\documentclass[letterpaper, preprint, paper,11pt]{AAS}	
\usepackage[colorlinks=true, pdfstartview=FitV, linkcolor=black, citecolor= black, urlcolor= black]{hyperref}
\usepackage{overcite}
\usepackage{footnpag}			      	

\usepackage{amsmath}
\usepackage{amsfonts}
\usepackage{amssymb}
\usepackage{siunitx}
\usepackage{booktabs}
\usepackage{mathtools}
\usepackage{caption}
\usepackage{subcaption}
\usepackage{multirow}
\usepackage{lscape}
\usepackage{upgreek}
\usepackage{algorithm}
\usepackage{algpseudocode}
\usepackage{todonotes}

\newcommand{\Frame}{\mathcal{F}}

\newcommand{\cbold}{\boldsymbol{c}}
\newcommand{\Ebold}{\boldsymbol{E}}
\newcommand{\Fbold}{\boldsymbol{F}}
\newcommand{\Gbold}{\boldsymbol{G}}

\newcommand{\Tbold}{\boldsymbol{T}}

\newcommand{\abold}{\boldsymbol{a}}

\newcommand{\rbold}{\boldsymbol{r}}
\newcommand{\vbold}{\boldsymbol{v}}
\newcommand{\ubold}{\boldsymbol{u}}
\newcommand{\Ubold}{\boldsymbol{U}}

\newcommand{\xbold}{\boldsymbol{x}}
\newcommand{\Xbold}{\boldsymbol{X}}

\newcommand{\fbold}{\boldsymbol{f}}
\newcommand{\Ibold}{\boldsymbol{I}}

\newcommand{\xibold}{\boldsymbol{\xi}}
\newcommand{\lambdabold}{\boldsymbol{\lambda}}
\newcommand{\Phibold}{\boldsymbol{\Phi}}

\newcommand{\minimize}[1]{\underset{#1}{\operatorname{min}}}

\newcommand*{\permcomb}[4][0mu]{{{}^{#3}\mkern#1#2_{#4}}}

\newcommand*{\comb}[1][-1mu]{\permcomb[#1]{C}}

\PaperNumber{25-717}

\begin{document}

\title{Multi-Vehicle Guidance for Formation Flight on Libration Point Orbits}

\author{Yuri Shimane\thanks{PhD Candidate, School of Aerospace Engineering, Georgia Institute of Technology, GA 30332, USA.},
Purnanand Elango\thanks{Research Scientist, Mitsubishi Electric Research Laboratories, Cambridge, MA 02139, USA.},
\ and Avishai Weiss\thanks{Senior Principal Research Scientist, Mitsubishi Electric Research Laboratories, Cambridge, MA 02139, USA.}
}

\maketitle{}

\begin{abstract}
The multiple spacecraft guidance problem for proximity flight in libration point orbit is considered.
A nonlinear optimal control problem with continuous-time path constraints enforcing minimum separation between each spacecraft is formulated.
The path constraints are enforced via an isoperimetric reformulation, and the problem is solved via a sequential convex programming.
The proposed approach does not necessitate specific dynamic system structures to provide continuous-time guarantees for minimum separation within a fuel-optimal solution.
The optimal control problem is deployed within a model predictive control scheme and demonstrated in the ephemeris model dynamics.
\end{abstract}

\section{Introduction}
Safe guidance schemes for formation flight around libration point orbits (LPOs) are of interest due to their relevance in constellation missions~\cite{Sherman2023,Capannolo2023} as well as for proximity operations and rendezvous with assets such as the Gateway along its near rectilinear halo orbit (NRHO)~\cite{Bucci2018,Bucchioni2020,Blazquez2020}.
Motivated by the success of the relative orbital elements-based formulation for formation flight in Earth orbits, the use of the center manifold has been studied for formation flying applications around LPOs.
The center manifold, which corresponds to a non-expanding, oscillatory eigenvector, provides a natural mechanism to remain in the vicinity of an LPO and may accommodate multiple vehicles by offsetting the phase of the oscillation.
For example, Calico and Wiesel~\cite{Calico1984} demonstrated the use of Floquet modes to provide insight into relative motion in the vicinity of LPOs. 
Linearized relative orbital elements were introduced by Hsiao and Scheeres~\cite{Hsiao2002} to arrive at a feedback control law that produced oscillatory motions that are similar to the motion in the center eigenspace, as studied in Scheeres et al.~\cite{Scheeres2003}. 
Recently, the local toroidal coordinate (LTC)~\cite{Elliott2022} has been gaining attention as a geometric framework for relative motion around LPOs in simplified dynamics models, such as the circular restricted three-body problem (CR3BP).
An LTC is defined using the non-expanding, oscillatory mode of the LPO's monodromy matrix and is well-suited to study bounded motion and guidance problems in the vicinity of the LPO. 
Elliott and Bosanac~\cite{Elliott2021} developed a targeting-based guidance scheme in the LTC defined in the CR3BP, which is also tested in a high-fidelity ephemeris model (HFEM).
Takubo et al.~\cite{Takubo2025} also adopted the LTC system and developed optimization-based guidance problems in the CR3BP that enforce passive safety through geometric constraints.

While QPT also exists in HFEM~\cite{Howell2005}, the lack of an exact center manifold in the HFEM prohibits the precise definition of LTC and, thereby, the direct construction of QPT.
Instead, a QPT from simplified models such as the CR3BP may be transitioned through differential correction~\cite{Howell2005}.
To ensure the transitioned QPT exhibits the desirable properties from the QPT in the simplified model, the transition process must be carefully monitored and tuned by the trajectory designer.
In the presence of uncertainties and corrective control maneuvers, as is the case in actual flight, the performance of QPT-based approaches may be empirically shown~\cite{Elliott2021,Takubo2025}.
Still, the inherent need for transitioning from simpler dynamics models renders any formal guarantees difficult to establish in the HFEM. 

As an alternative to LTC-based guidance, we formulate a general multi-vehicle guidance problem (MVGP) that is agnostic to the existence of a center manifold.
The MVGP is a nonlinear model predictive control (MPC) scheme, which extends a recently developed MPC-based station-keeping scheme for a single spacecraft following an LPO~\cite{Shimane2025} to the multi-vehicle case.
As in the single spacecraft case, we adopt an impulsive control model with subsequent maneuvers placed one revolution apart, abiding by typical operational requirements on Earth-Moon LPOs.
Desired properties, such as bounded inter-spacecraft range between member spacecraft in the formation, are incorporated into the MVGP as \textit{continuous-time} (CT) path constraints, enforced at all times within the MPC's prediction horizon.
The propellant-optimal solution to the MVGP will naturally favor toroidal configurations, where the spacecraft remains in the vicinity of the LPO while ensuring a minimum separation between themselves.
Even in the HFEM with no exact central manifold and in the presence of uncertainty, where a toroidal configuration is not strictly non-expanding and may be disturbed, the MVGP will yield the propellant-minimizing solution that ensures a minimum separation while keeping all spacecraft in the vicinity of the LPO.
In general, the MVGP is agnostic to natural dynamical structures, thereby allowing the incorporation of any arbitrary constraints necessitated by operations.

In order to maintain the constellation over an extended duration, the MPC is re-instantiated and solved at each revolution along the LPO, sliding the prediction horizon downstream by one revolution.
In this recursive process, due to state estimation, control execution, and dynamics modeling error, a feasible solution to the MVGP where the path constraint is active may arrive at a realization that lies closer than the required minimum separation at the subsequent control time, thus resulting in a re-instantiated MVGP that has no feasible solution.
To circumvent this issue and ensure recursive feasibility under uncertainty, we enforce separation thresholds that monotonically increase across the MPC's control horizon.

To handle the CT path constraints, which are, by nature, functional inequality constraints, we make use of their isoperimetric form to enforce these by augmenting the dynamics~\cite{Teo1987}.
Upon reformulation of the CT path constraints, the MVGP is a non-convex nonlinear program (NLP), which we solve via sequential convex programming (SCP).
Discussions on the property of incorporating isoperimetric constraints into the NLP are treated in detail in Elango et al.~\cite{Elango2024}. 
The incorporation of isoperimetric constraints with SCP has been applied to the powered descent and landing problem~\cite{Elango2025}, cislunar spacecraft rendezvous~\cite{Elango2025RDV}, and for geostationary satellite station-keeping~\cite{Pavlasek2025}.
Notably, Pavlasek et al.~\cite{Pavlasek2025} considers an MPC scheme where the SCP is solved within a recursive horizon, similarly to this work.
In lieu of the prox-linear method~\cite{Lewis2008,Drusvyatskiy2018} adopted in~\cite{Elango2024,Elango2025}, which involves an exact penalization of non-convex constraints, we use the Augmented-Lagrangian-based SCP algorithm~\cite{Oguri2023}, which consists of updating both optimization variables and Lagrange multipliers through a primal-dual formalism.
Both SCP algorithms possess theoretical convergence guarantees to a local feasible minimizer of the original non-convex NLP, and are thus reasonable choices for solving the MVGP; in the present work, the augmented Lagrangian SCP is chosen due to its successful application to trajectory design problems in cislunar space~\cite{Kumagai2024,Kumagai2024AAS}.

\section{Background}
We first provide background on the dynamics, then briefly introduce the notion of the baseline NRHO.

\subsection{Dynamics}
We consider the N-body HFEM dynamics for the spacecraft.
We use an inertial frame~$\Frame_{\rm Inr}$ that corresponds to the J2000 frame, centered on the Moon, where we resolve the dynamics of the spacecraft.
In addition, we introduce the Earth-Moon rotating frame~$\Frame_{\rm EM}$ also centered on the Moon, with its first axis aligned with the instantaneous Earth-to-Moon direction, and its second axis aligned with the velocity direction of the Moon with respect to the Earth. 

Let $\rbold \in \mathbb{R}^3$ denote the position vector in~$\Frame_{\rm Inr}$, and $\vbold \triangleq \dot{\rbold} \in \mathbb{R}^3$ denote the velocity in~$\Frame_{\rm Inr}$.
Let $\xbold = [\rbold^T, \vbold^T]^T \in \mathbb{R}^6$ represent the spacecraft state.
The natural dynamics is given by
\begin{equation}    \label{eq:dynamics}
    \dot{\xbold}(t) = \begin{bmatrix}
        \dot{\rbold}(t) \\ \dot{\vbold}(t)
    \end{bmatrix} = 
    \begin{bmatrix}
        \vbold(t) \\
        -\dfrac{\mu}{\|\rbold(t)\|_2^3}\rbold(t)
        + \sum_b \abold_b(t)
        + \abold_{\rm SRP}(t)
    \end{bmatrix}
    ,
\end{equation}
where $\mu$ is the Moon's gravitational parameter, $\abold_b(t)$ is the third-body acceleration of body $b$ with gravitational parameter $\mu_b$ and position vector $\rbold_b$ with respect to the Moon, 
\begin{equation}    
    \abold_b(t) = - \mu_b \left( \dfrac{\rbold(t) - \rbold_b(t)}{\|\rbold(t) - \rbold_b(t)\|_2^3} + \dfrac{\rbold_b(t)}{\|\rbold_b(t)\|_2^3} \right)
    ,
\end{equation}
where we consider third-body perturbations due to the Earth and the Sun,
and $\abold_{\rm SRP}$ is the acceleration due to solar radiation pressure (SRP), where the cannonball model is assumed~\cite{Vallado2001},
\begin{equation}
    \abold_{\rm SRP} = 
    P_{\rm \odot} \left( \dfrac{\mathrm{AU}}{\|\rbold - \rbold_{\odot}\|_2} \right)^2 \dfrac{C_r A}{m} \dfrac{\rbold - \rbold_{\odot}}{\|\rbold - \rbold_{\odot}\|_2}
    ,
\end{equation}
where $P_{\rm \odot}$ is the Sun's pressure at $1$ astronomical unit $\mathrm{AU} = 149.6\times10^6$ \SI{}{km}, $C_r$ is the radiation pressure coefficient, $A/m$ is the pressure area-to-mass ratio of the spacecraft, and $\rbold_{\odot}$ is the position vector of the Sun with respect to the Moon.

In this work, we assume impulsive control $\ubold(t_{k|j}) \in \mathbb{R}^3$, imparting an instantaneous change in $\vbold$ at time $t_{k|j}$.
The spacecraft state at time $t$ is thus obtained with the integration
\begin{equation}    \label{eq:IVP}
    \xbold(t) = \xbold(t_{k|j})
    + 
    \begin{bmatrix}
        \boldsymbol{0}_{3\times3} \\ \Ibold_3
    \end{bmatrix} \ubold(t_{k|j})
    + \int_{t_{k|j}}^t \fbold(\xbold(\tau),\tau) \mathrm{d}\tau
    ,
\end{equation}
where $\delta$ is the Dirac delta function.
The Jacobian of~\eqref{eq:dynamics} is given by
\begin{equation}    \label{eq:Jacobian}
    \dfrac{\partial \fbold}{\partial \xbold} 
    = \begin{bmatrix}
        \boldsymbol{0}_{3\times3} & \Ibold_3 \\
        \Gbold & \boldsymbol{0}_{3\times3}
    \end{bmatrix}
    ,
\end{equation}
where $\Gbold$ is the gravity gradient matrix for the central body and third-body perturbations, and partials with respect to the solar radiation pressure.

\subsubsection{Canonical Scales}
To facilitate solving NLPs involving $\xbold$ and $\ubold$, we rescale the dynamics into canonical scales. Based on a user-defined distance unit $\rm DU$, we define the velocity unit $\rm VU \triangleq \sqrt{\mu/\mathrm{DU}}$ and time unit $\rm TU \triangleq {\mathrm{DU}/\mathrm{VU}}$,
such that $\mu$ of the central body in~\eqref{eq:dynamics} is rescaled to $1$.

\subsubsection{Multi-Vehicle Augmented System}
In this work, we consider formation flight scenarios involving $M \geq 2$ spacecraft flying in the vicinity of a baseline LPO.
Let $\Xbold \in \mathbb{R}^{6M}$ and $\Ubold \in \mathbb{R}^{3 M}$ represent the concatenated states and controls of $M$ spacecraft,
\begin{equation}    \label{eq:augmented_state_control}
    \Xbold(t) = \begin{bmatrix}
        \xbold_{0}(t) \\ \vdots \\ \xbold_{M-1}(t)
    \end{bmatrix}
    ,\quad
    \Ubold(t) = \begin{bmatrix}
        \ubold_{0}(t) \\ \vdots \\ \ubold_{M-1}(t)
    \end{bmatrix}
    ,
\end{equation}
where $\xbold_i$ and $\ubold_i$ correspond to the state and control of the $i^{\rm th}$ spacecraft, where $i\in\mathcal{I} = \{0,\ldots,M-1\}$.
The natural dynamics for $\Xbold$ is given by
\begin{equation}    
    \dot{\Xbold}(t) = \Fbold(\Xbold(t),t) =
    \begin{bmatrix}
        \fbold_{0}(\xbold_{0}(t),t) \\ \vdots \\ \fbold_{M-1}(\xbold_{M-1}(t),t)
    \end{bmatrix}
    ,
\end{equation}
where $\fbold_i$ is the dynamics~\eqref{eq:dynamics} for the $i^{\rm th}$ spacecraft, with corresponding Jacobian
\begin{equation}    
    \dfrac{\partial \Fbold}{\partial \Xbold}
    =
    \operatorname{block-diag}\left(
        \dfrac{\partial \fbold_{0}}{\partial \xbold_{0}},
        \cdots,
        \dfrac{\partial \fbold_{M-1}}{\partial \xbold_{M-1}}
    \right)
    ,
\end{equation}
where $\partial \fbold_i/\partial \xbold_i$ is the Jacobian~\eqref{eq:Jacobian} about the dynamics of the $i^{\rm th}$ spacecraft.

\subsection{Baseline NRHO}
We study the formation flight guidance problem relative to a predefined \textit{baseline}, a precomputed deterministic trajectory in the HFEM, which the formation is expected to recursively track.
While in general, each member spacecraft within the formation may have unique baselines, we consider in this work that the same baseline is tracked by all member spacecraft.
We specifically use the 15-year-long NRHO by NASA~\cite{Lee2019} as the baseline.
This NRHO, originally designed as a preliminary orbit for the Gateway, has an orbital period of approximately $6.55$ \SI{}{days}, exhibiting approximately a 9:2 resonance with the Sun-Earth-Moon synodic month.
In subsequent sections, let $(\cdot)_{\rm bsln}$ denote a quantity along the baseline.


\section{Multi-Vehicle Guidance Problems}
\label{sec:MVGP}
In this Section, we develop the MVGP formulation.
We begin with a brief discussion on the targeting MPC approach to station-keeping, developed by the authors for the case with a single spacecraft~\cite{Shimane2025}.
We then define the discretization, multiple shooting, and targeting constraint with the augmented state for the multi-vehicle case.
Next, we formulate the MVGP as a discretized non-convex optimization with continuous-time path constraints, which take the form of functional inequality constraints.
Finally, we consider two distinctions in which we consider MVGP instances.
The first distinction is between rendezvous vs. formation flying constellation applications, which translates to differences in the initial states and tracked references of each spacecraft.
The second distinction is between cooperative vs. hierarchical scenarios: in a cooperative scenario, control sequences for all vehicles are simultaneously computed, while in a hierarchical scenario, control sequence(s) of some vehicle(s) are fixed to values obtained a priori, and the control sequence(s) of remaining vehicle(s) subject to constraints in the MVGP are computed.

\subsection{Station-Keeping as a Model Predictive Control Problem}
Due to the existence of unstable modes, spacecraft operating on LPOs must conduct station-keeping to remain in the vicinity of a precomputed reference trajectory, or \textit{baseline}, despite the presence of uncertainties.
On an LPO with an orbital period of a few days, a typical operation requirement is to conduct at most one station-keeping maneuver per revolution.
The Gateway, to be operated on an NRHO with an orbital period of $\approx 6.55$ \SI{}{days}, is subjected to this requirement.
Whereas several state-of-the-art station-keeping algorithms involve designing the single maneuver at a given revolution to target some condition downstream, Shimane et al.~\cite{Shimane2025} introduce an MPC approach, where a control horizon incorporates not only the immediate but also subsequent maneuvers.
The extension of the control horizon recovers the controllability for full-state targeting, thus allowing for the spacecraft to recursively track the baseline both in the geometry of the NRHO and in phase.
The MVGP in this work is an extension of the MPC scheme to the multi-vehicle case.

\subsubsection{Targeting Constraint}
Let $(\cdot)_{k|j}$ denote the quantity predicted at $k \in \mathcal{K} = \{ 0,\ldots,N-1 \}$ time increments ahead of the current state estimate $\xbold(t_j)$ at time $t_j$.
We consider a prediction horizon consisting of $N - 1$ revolutions, with $N$ control nodes spaced one revolution apart at times $t_{0|j},\ldots,t_{N-1|j}$, where $\ubold_{0|j},\ldots,\ubold_{N-1|j}$ denote the corresponding maneuvers.
The MPC seeks to target the baseline at $t_{N-1|j}$,
\begin{equation}
    \xbold(t_{N-1|j}) \in \mathcal{X}(t_{N-1|j}),
\end{equation}
while minimizing the cumulative control effort.
One simple yet effective choice of $\mathcal{X}(t_{N-1|j})$ is an ellipsoid centered at the baseline state at $t_{\rm final}$, denoted by $\xbold_{\rm bsln}(t_{\rm final})$, with radii $\epsilon_r$ in position and $\epsilon_v$ in velocity~\cite{Shimane2025}.
Then, the targeting constraint becomes a pair of second-order cone constraints
\begin{align}   \label{eq:generic_targeting_constraint}
    \| \rbold(t_{N-1|j}) - \rbold_{\rm bsln}(t_{N-1|j}) \|_2 \leq \epsilon_r
    ,\quad
    \| \vbold(t_{N-1|j}) - \vbold_{\rm bsln}(t_{N-1|j}) \|_2 \leq \epsilon_v
    .
\end{align}
Once the MPC is solved and maneuver $\ubold_{0|j}$ is executed\footnote{In reality, a command to execute $\ubold_{0|j}$ is sent, but an imperfect maneuver $\ubold_{0|j} + \delta \ubold$ is imparted to the spacecraft.}, the MPC is re-invoked at time $t_{j+1}$, shifting the horizon by a revolution. 
The size of the targeting ellipsoid in~\eqref{eq:generic_targeting_constraint}, $\epsilon_r$ and $\epsilon_v$ are tuning parameters station-keeping scheme, whose appropriate value that yields the fuel-optimal control depends on the level of uncertainties and length of prediction horizon $N$.
For more discussion on the choice of $\epsilon_r$ and $\epsilon_v$, see Shimane et al.~\cite{Shimane2025}.

The control times $t_{0|j},\ldots,t_{N-1|j}$ are fixed based on selecting an osculating true anomaly where the maneuvers are to occur, $\theta_{\rm man}$.
We then set $t_{k|j}$ such that the osculating true anomaly of the baseline, $\theta_{\rm bsln}$, defined as 
\begin{equation}    
\begin{aligned}
    &\theta_{\rm bsln} = \operatorname{atan2}\left(
        h_{\rm bsln} v_{r,\rm bsln}, h_{\rm bsln}^2/ \|\rbold_{\rm bsln}\|_2 - \mu
    \right),
    \\&
    h_{\rm bsln} = \| \rbold_{\rm bsln} \times \vbold_{\rm bsln} \|_2,
    \quad
    v_{r,\rm bsln} = \rbold_{\rm bsln}^T \vbold_{\rm bsln} / \| \rbold_{\rm bsln} \|_2,
\end{aligned}
\end{equation}
is equal to $\theta_{\rm man}$ along the $k^{\rm th}$ revolution.

\subsection{Problem Formulation}
Let $(\cdot)_{i,k|j}$ denote a quantity associated with the $i^{\rm th}$ spacecraft at $t_{k|j}$.
Following the notations from~\eqref{eq:augmented_state_control}, the augmented states and controls at $t_{k|j}$ are denoted by
\begin{equation}    
    \Xbold_{k|j} = \begin{bmatrix}
        \xbold_{0,k|j} \\ \vdots \\ \xbold_{M-1,k|j}
    \end{bmatrix}
    ,\quad
    \Ubold_{k|j} = \begin{bmatrix}
        \ubold_{0,k|j} \\ \vdots \\ \ubold_{M-1,k|j}
    \end{bmatrix}
    .
\end{equation}
Let $\mathcal{C}$ denote the set of pairs $(i,j)$ where $i,j \in \mathcal{I}$ and $i \neq j$.
We extend the single vehicle guidance problem from Shimane et al.~\cite{Shimane2025} by steering the concatenated states of $M$ spacecraft with corresponding $M$ controls from~\eqref{eq:augmented_state_control}, also incorporating path constraints on the inter-satellite ranges.
The resulting MVGP is
\begin{subequations}    \label{eq:MVGP}
\begin{align}
    \minimize{
        \substack{\Xbold_{0|j},\ldots,\Xbold_{N-1|j}\\\Ubold_{0|j},\ldots,\Ubold_{N-1|j}}
    }
    \quad& \sum_{i=0}^{M-1} \sum_{k=0}^{N-1} \| \ubold_{i,k} \|_2
        \label{eq:MVGP_obj}
    \\\text{s.t.}\quad&
    \Xbold_{k+1|j} = \Xbold_{k|j} + \Ebold \Ubold_{k|j} 
    + \int_{t_{k|j}}^{t_{k+1|j}} \Fbold(\Xbold(t),t)\mathrm{d}t
        \quad \forall k \in \mathcal{K} \setminus \{N-1\} 
    \label{eq:MVGP_dynamics}
    \\&
    \Xbold_{0|j} = \hat{\Xbold}(t_j)
    \label{eq:MVGP_initialconditions}
    \\&
    \| \rbold_{i,N-1|j} - \rbold_{\mathrm{bsln},i}(t_{N-1|j}) \|_2 \leq \epsilon_r
        \quad \forall i \in \mathcal{I} 
    \label{eq:MVGP_finalset_r}
    \\&
    \| \vbold_{i,N-1|j} + \ubold_{i,N-1|j} - \vbold_{\mathrm{bsln},i}(t_{N-1|j}) \|_2 \leq \epsilon_v
        \quad \forall i \in \mathcal{I} 
    \label{eq:MVGP_finalset_v}
    \\&
    \| \rbold_{i}(t) - \rbold_{j}(t) \|_2 \geq \Delta r_{\min}
        \quad \forall
        (i,j) \in \mathcal{C},\,
        t\in[t_{0|j},t_{N-1|j}]
    \label{eq:MVGP_minsep}
    \\&
    \| \rbold_{i}(t) - \rbold_{j}(t) \|_2 \leq \Delta r_{\max}
        \quad \forall 
        (i,j) \in \mathcal{C},\,
        t\in[t_{0|j},t_{N-1|j}]
    \label{eq:MVGP_maxsep}
\end{align}
\end{subequations}
where in~\eqref{eq:MVGP_dynamics}, $\Ebold$ is given by
\begin{equation}
    \Ebold = \operatorname{block-diag}
    \underbrace{
        \left(
            \begin{bmatrix} \boldsymbol{0}_{3 \times 3} \\ \Ibold_3 \end{bmatrix}, \ldots, 
            \begin{bmatrix} \boldsymbol{0}_{3 \times 3} \\ \Ibold_3 \end{bmatrix}
        \right)
    }_{\text{$M$ times}}
    .
\end{equation}
The objective~\eqref{eq:MVGP_obj} minimizes the sum of the maneuver costs of all spacecraft.
Constraints~\eqref{eq:MVGP_dynamics} ensure continuity of the dynamics;
constraint~\eqref{eq:MVGP_initialconditions} enforces the initial state to coincide with the state estimate of each spacecraft at $t_j$, concatenated as $\hat{\Xbold}(t_j)$; constraints~\eqref{eq:MVGP_finalset_r} and~\eqref{eq:MVGP_finalset_v} ensure each spacecraft arrives in the ellipsoidal target set $\mathcal{X}(t_{N-1|j})$ with respect to the baseline orbit for the $i^{\rm th}$ spacecraft; finally, path constraints~\eqref{eq:MVGP_minsep} and~\eqref{eq:MVGP_maxsep} ensure the inter-spacecraft distances are bounded from below and above by $\Delta r_{\min}$ and $\Delta r_{\max}$.

\section{Solution Approach}
\label{sec:solution_approach}
In this Section, we present the solution approach to solve problem~\eqref{eq:MVGP}.
The approach consists of two steps: first, we demonstrate the use of isoperimetric forms~\cite{Teo1987} for the continuous-time path constraints~\eqref{eq:MVGP_minsep} to recast~\eqref{eq:MVGP} as a non-convex NLP with no functional constraints.
Then, we provide a summary of the augmented Lagrangian-based SCP algorithm~\cite{Oguri2023}, which is used to solve the non-convex NLP by solving a sequence of convex programs.
Finally, we provide discussions on scaling and tuning hyperparameters pertaining to the isoperimetric reformulation and the SCP algorithm.

\subsection{Gradual constraint tightening for Recursively Feasible Solutions}
The MVGP~\eqref{eq:MVGP} is recursively solved within an MPC scheme.
Identically to the single-vehicle station-keeping MPC, at a given time $t_j$, the earliest controls $\Ubold_{0|j}$, corrupted with noise, are imparted on the spacecraft. The states are then propagated until $t_{j+1}$, at which point the MPC is re-instantiated to obtain the next control to be executed, $\Ubold_{0|j+1}$.
The true state at $t_{j+1}$, $\Xbold(t_{j+1})$, following the corrupted controls $\Tilde{\Ubold}_{0|j}$ at $t_j$, is
\begin{equation}    \label{eq:recurse_MPC_onestep}
    \Xbold(t_{j+1})
    =
    \Xbold(t_j) + \Ebold \Tilde{\Ubold}_{0|j} + \int_{t_j}^{t_{j+1}} \Fbold(\Xbold(t),t)\mathrm{d}t
    ,
\end{equation}
and the estimated state at $t_{j+1}$ is
\begin{equation}   \label{eq:state_estimate}
    \hat{\Xbold}(t_{j+1}) = \Xbold(t_j) + \delta \Xbold_{\rm nav}
\end{equation}
where $\delta \Xbold_{\rm nav} \in \mathbb{R}^6$ is some state estimation error.
During this recursion, even if $\Ubold_{0|j}$ was a feasible solution of~\eqref{eq:MVGP} at $t_j$, at $t_{j+1}$, the initial conditions $\hat{\Xbold}(t_{j+1})$ enforced by~\eqref{eq:MVGP_initialconditions} may result in path constraints~\eqref{eq:MVGP_minsep} and~\eqref{eq:MVGP_maxsep} that cannot be satisfied.

\begin{figure}[t]
    \centering
    \includegraphics[width=0.95\linewidth]{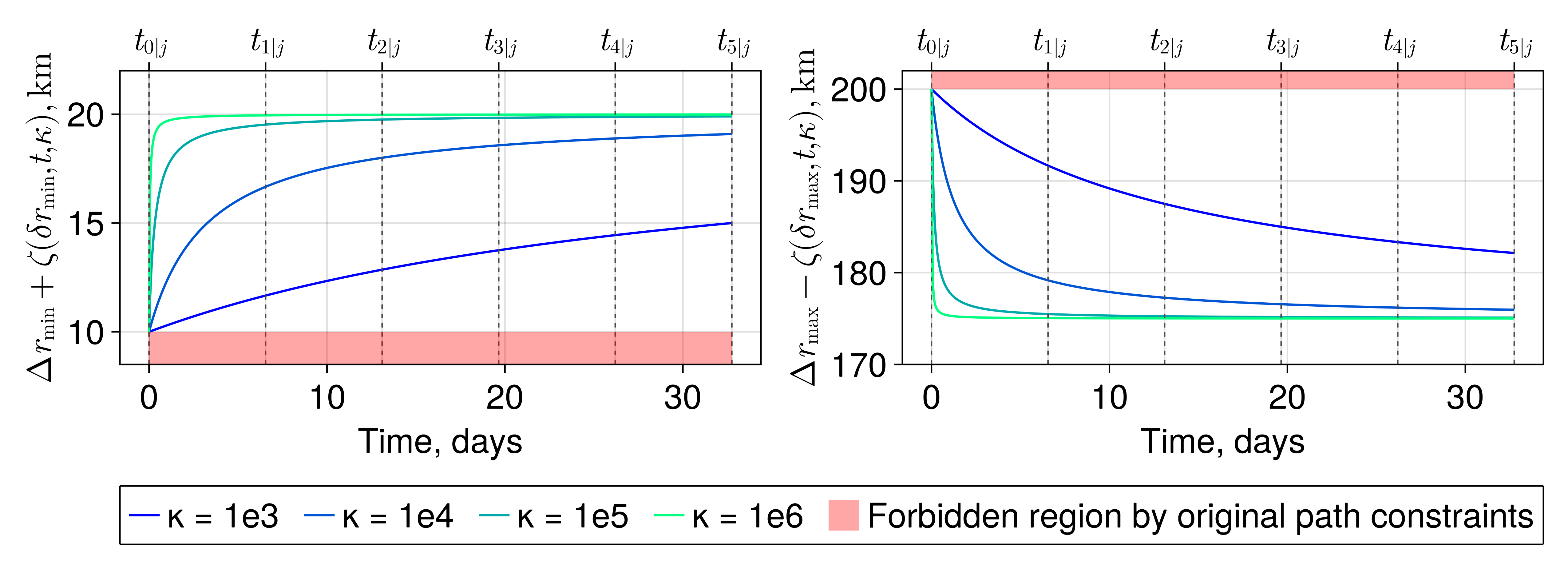}
    \caption{Time-dependent bounds on inter-spacecraft separation path constraints for various tuning parameter $\kappa$ with $\Delta r_{\min} = 10$ \SI{}{km}, $\delta r_{\min} = 10$ \SI{}{km}, $\Delta r_{\max} = 200$ \SI{}{km}, $\delta r_{\max} = 25$ \SI{}{km}}
    \label{fig:range_bounds_tightened}
\end{figure}

One approach would be to reformulate the MVGP~\eqref{eq:MVGP_maxsep} using chance constraints, incorporating uncertainties due to state estimation error, model mismatch, and control execution error.
Chance constraints reformulated in deterministic form have been considered in several past works on spacecraft guidance problems~\cite{Kumagai2024AAS,Lew2020,Takubo2024-tf,Oguri2024-of,Elango2025RDV}.
In this work, we instead opt for an approximate approach, where the right-hand sides of the path constraints are tightened along the prediction horizon of the MPC.
The minimum separation constraint~\eqref{eq:MVGP_minsep} is modified to
\begin{equation}    \label{eq:minsep_tightened}
    \| \rbold_i(t) - \rbold_j(t) \|_2 \geq \Delta r_{\min} + \zeta(\delta r_{\min}, t, \kappa)
    ,
\end{equation}
and the maximum separation constraint~\eqref{eq:MVGP_maxsep} is modified to
\begin{equation}    \label{eq:maxsep_tightened}
    \| \rbold_i(t) - \rbold_j(t) \|_2 \leq \Delta r_{\max} - \zeta(\delta r_{\max}, t, \kappa)
    ,
\end{equation}
where the function $\zeta: \mathbb{R}^{3} \to \mathbb{R}$ is given by
\begin{equation}    \label{eq:zeta_definition}
    \zeta (\delta r, t, \kappa) = \delta r - \dfrac{1}{(\kappa \tau(t) + 1 / \delta r)}
    ,\quad
    \tau(t) = \dfrac{t - t_{0|j}}{t_{N-1|j} - t_{0|j}}
    ,
\end{equation}
where $\delta r_{\min}$ and $\delta r_{\max}$ are additional separation thresholds, and $\kappa$ is a tuning parameter.
Figure~\ref{fig:range_bounds_tightened} shows traces where the tightened constraints~\eqref{eq:minsep_tightened} and~\eqref{eq:maxsep_tightened} are active for $\kappa = 10^3 \sim 10^6$.
In order to ensure recursive feasibility, $\kappa$, $\delta r_{\min}$, and $\delta r_{\max}$ must be chosen such that for any error realization, $\hat{\Xbold}(t_{j+1})$ by~\eqref{eq:state_estimate} lies outside of the forbidden region at the next iteration $t_{j+1}$, which coincides with $t_{1|j}$.
These parameters are empirically tuned in a preliminary experiment, and values are reported in the Numerical Results section.

\subsection{Continuous-Time Path Constraints}
The functional inequality constraints~\eqref{eq:MVGP_minsep} and~\eqref{eq:MVGP_maxsep} cannot be enforced in their current form.
One approximation is to enforce the constraints at the nodes,
\begin{subequations}
\begin{align}
    \| \rbold_{i,k|j} - \rbold_{j,k|j} \|_2 \geq \Delta r_{\rm min} + \zeta(\delta r_{\min}, t_{k|j}, \kappa)
        \quad \forall (i,j) \in \mathcal{C},\,
        k \in \mathcal{K} 
    ,   \label{eq:minsep_nodes}
    \\ 
    \| \rbold_{i,k|j} - \rbold_{j,k|j} \|_2 \leq \Delta r_{\rm max} - \zeta(\delta r_{\max}, t_{k|j}, \kappa)
        \quad \forall (i,j) \in \mathcal{C},\,
        k \in \mathcal{K} 
    ,   \label{eq:maxsep_nodes}
\end{align}
\end{subequations}
but at the risk of inter-sample constraint violations.
Instead, we recast~\eqref{eq:MVGP_minsep} and~\eqref{eq:MVGP_maxsep} in isoperimetric form, which can be incorporated as part of an entirely finite-dimensional NLP.
We may rewrite a generic inequality path constraint $g$ to its isoperimetric form,
\begin{equation}    \label{eq:ct_constraints_isoperimetric_equivalence}
    g(\Xbold(t),\Ubold(t),t) \leq 0
    \quad t \in [t_{k|j},t_{k+1|j}]
    \quad \Leftrightarrow \quad 
    \int_{t_{k|j}}^{t_{k+1|j}}
    \left| g(\Xbold(t),\Ubold(t),t) \right|^{\alpha}_+ \mathrm{d}t = 0
    .
\end{equation}
where $\alpha > 1$.
The reformulation of continuous-time path constraints via~\eqref{eq:ct_constraints_isoperimetric_equivalence} is analyzed in detail in Elango et al.~\cite{Elango2024}, and has been applied for three- and six-degrees-of-freedom rocket landing problems~\cite{Elango2024,Chari2024,Elango2025}, as well as for station-keeping on GEO~\cite{Pavlasek2025}.

In the case of the MVGP, $g$ corresponds to the separation constraints~\eqref{eq:MVGP_minsep} and~\eqref{eq:MVGP_maxsep} for each spacecraft $(i,j)$, given by
\begin{equation}
    \label{eq:pc_MVGP_integral}
    \int_{t_{k|j}}^{t_{k+1|j}}
    \left| g^{(i,j)}_{\Delta r}(\Xbold(t),t) \right|^{\alpha}_+ \mathrm{d}t = 0
    ,
\end{equation}
where
\begin{equation}   \label{eq:integrand_minsep}
    g^{(i,j)}_{\Delta r_{\min}}(\Xbold(t),t) = \Delta r_{\min} + \zeta(\delta r_{\min},t,\kappa) - \| \rbold_i(t) - \rbold_j(t) \|_2
    ,
\end{equation}
for~\eqref{eq:MVGP_minsep} and
\begin{equation}    \label{eq:integrand_maxsep}
    g^{(i,j)}_{\Delta r_{\max}}(\Xbold(t),t) = -\Delta r_{\max} + \zeta(\delta r_{\max},t,\kappa) + \| \rbold_i(t) - \rbold_j(t) \|_2
    ,
\end{equation}
for~\eqref{eq:MVGP_maxsep}.
Note that there are $q = \comb{M}{2}$ constraint~\eqref{eq:integrand_minsep} and~\eqref{eq:integrand_maxsep}, respectively.
In order to enforce~\eqref{eq:pc_MVGP_integral}, we introduce slack variables $y^{(i,j)}_{\Delta r_{\min}}$ and $y^{(i,j)}_{\Delta r_{\max}}$ for each constraint, and define the augmented state $\Xbold^{\rm g} \in \mathbb{R}^{6M + 2q}$
\begin{equation}    
    \Xbold^{\rm g}(t) = \begin{bmatrix}
        \Xbold(t) \\
        y^{(1,2)}_{\Delta r_{\min}}(t) \\ \vdots \\
        y^{(M-1,M)}_{\Delta r_{\min}}(t) 
        \\[0.4em]
        y^{(1,2)}_{\Delta r_{\max}}(t) \\ \vdots \\
        y^{(M-1,M)}_{\Delta r_{\max}}(t)
    \end{bmatrix}
    .
\end{equation}
The dynamics for $\Xbold^{\rm g}$ is given by
\begin{equation}    
    \dot{\Xbold}^{\rm g}(t) =
    \Fbold^{\rm g} (\Xbold^{\rm g},t) =
    \begin{bmatrix}
        \Fbold(\Xbold,t) \\
        \Gbold_{\Delta r_{\min}}(\Xbold,t) \\
        \Gbold_{\Delta r_{\max}}(\Xbold,t)
    \end{bmatrix}
    ,\quad
    \Gbold_{\Delta r}(\Xbold,t) = 
    \begin{bmatrix}
        \left| g^{(1,2)}_{\Delta r}(\Xbold(t),t) \right|^{\alpha}_+ \\ 
        \vdots \\
        \left| g^{(M-1,M)}_{\Delta r}(\Xbold(t),t) \right|^{\alpha}_+
    \end{bmatrix}
    ,
\end{equation}
where $(\cdot)_{\Delta r}$ denote quantities associated with either $g_{\Delta r_{\min}}$ or $g_{\Delta r_{\max}}$, respectively.
The Jacobian for the augmented system is given by
\begin{equation}    
    \dfrac{\partial \Fbold^{\rm g}}{\partial \Xbold^{\rm g}}
    = 
    \begin{bmatrix}
        \dfrac{\partial \Fbold}{\partial \Xbold} 
        & 
        \boldsymbol{0}_{6M \times q}
        \\[0.6em]
        \dfrac{\partial \Gbold_{\Delta r_{\min}}}{\partial \Xbold} & \boldsymbol{0}_{q \times q}
        \\[0.6em]
        \dfrac{\partial \Gbold_{\Delta r_{\max}}}{\partial \Xbold} & \boldsymbol{0}_{q \times q}
    \end{bmatrix}
    ,
    \quad
    \dfrac{\partial \Gbold_{\Delta r}}{\partial \Xbold} = 
    \begin{bmatrix}
        \dfrac{\partial \left|g^{(1,2)}_{\Delta r}\right|^{\alpha}_+ }{\partial \Xbold} & \boldsymbol{0}_{1 \times q}
        \\ \vdots \\
        \dfrac{\partial \left|g^{(M-1,M)}_{\Delta r}\right|^{\alpha}_+ }{\partial \Xbold} & \boldsymbol{0}_{1 \times q}
    \end{bmatrix}
    ,
\end{equation}
where the partials $\partial |g^{(i,j)}_{\Delta r_{\min}}|^{\alpha}_+ / \partial \Xbold$ and $\partial |g^{(i,j)}_{\Delta r_{\max}}|^{\alpha}_+ / \partial \Xbold$ is given by
\begin{equation}    
\begin{aligned}
    \dfrac{\partial |g^{(i,j)}_{\Delta r_{\min}}|^{\alpha}_+ }{\partial \Xbold}
    &=
    \alpha \left|
        \Delta r_{\min} + \zeta(t) - \| \rbold_i(t) - \rbold_j(t) \|_2
    \right|_+^{\alpha-1}
    \begin{bmatrix}
        G_{0|j} & \cdots & G_{M}
    \end{bmatrix}
    ,
    \\
    \dfrac{\partial |g^{(i,j)}_{\Delta r_{\max}}|^{\alpha}_+ }{\partial \Xbold}
    &=
    \alpha \left|
        -\Delta r_{\max} + \zeta(t) + \| \rbold_i(t) - \rbold_j(t) \|_2
    \right|_+^{\alpha-1}
    \begin{bmatrix}
        -G_{0|j} & \cdots & -G_{M}
    \end{bmatrix}
    ,
    \\
    G_{l} &= \begin{cases}
        \begin{bmatrix}
            \dfrac{-(\rbold_i(t) - \rbold_j(t))^T}{\| \rbold_i(t) - \rbold_j(t) \|_2}
            &
            \boldsymbol{0}_{1\times3}
        \end{bmatrix}
        & l = i,
        \\[0.6em]
        \begin{bmatrix}
            \dfrac{(\rbold_i(t) - \rbold_j(t))^T}{\| \rbold_i(t) - \rbold_j(t) \|_2}
            &
            \boldsymbol{0}_{1\times3}
        \end{bmatrix}
        & l = j,
        \\[0.3em]
        \boldsymbol{0}_{1\times6} & \text{otherwise}.
    \end{cases}
\end{aligned}
\end{equation}

Let $\Xbold^{\rm g}_{k|j}$ denote the augmented state on the $k^{\rm th}$ control node at $t_{k|j}$, with components $y^{(i,j)}_{\Delta r_{\min},k|j}$ and $y^{(i,j)}_{\Delta r_{\min},k|j}$ corresponding to the slacks $y^{(i,j)}_{\Delta r_{\min}}$ and $y^{(i,j)}_{\Delta r_{\max}}$ at $t_{k|j}$.
With the path constraints expressed in the isoperimetric form, the MVGP becomes
\begin{subequations}    \label{eq:MVGP_ct}
\begin{align}
    \minimize{
        \substack{\Xbold^{\rm g}_{0|j},\ldots,\Xbold^{\rm g}_{N-1|j}\\\Ubold_{0|j},\ldots,\Ubold_{N-1|j}}
    }
    \quad& \sum_{i=0}^{M-1} \sum_{k=0}^{N-1} \| \ubold_{i,k} \|_2
        \label{eq:MVGP_ct_obj}
    \\\text{s.t.}\quad&
    \Xbold^{\rm g}_{k+1|j} = \Xbold^{\rm g}_{k|j} 
    + \begin{bmatrix}
        \Ebold \\ \boldsymbol{0}_{q \times 3M}
    \end{bmatrix} \Ubold_{k|j}
    + \int_{t_{k|j}}^{t_{k+1|j}} \Fbold^{\rm g}(\Xbold^{\rm g}(t),t)\mathrm{d}t
        \quad \forall k \in \mathcal{K} \setminus \{N-1\} 
    \label{eq:MVGP_ct_dynamics}
    \\&
    \Xbold_{0|j} = \hat{\Xbold}(t_{0|j})
    \label{eq:MVGP_ct_initialconditions}
    \\&
    \| \rbold_{i,N-1|j} - \rbold_{\mathrm{bsln},i}(t_{N}) \|_2 \leq \epsilon_r
        \quad \forall i \in \mathcal{I} 
    \label{eq:MVGP_ct_finalset_r}
    \\&
    \| \vbold_{i,N-1|j} + \ubold_{i,N-1|j} - \vbold_{\mathrm{bsln},i}(t_{N}) \|_2 \leq \epsilon_v
        \quad \forall i \in \mathcal{I} 
    \label{eq:MVGP_ct_finalset_v}
    \\&
    y^{(i,j)}_{\Delta r_{\min},k|j} = y^{(i,j)}_{\Delta r_{\min},k+1|j}
        \quad \forall (i,j) \in \mathcal{C},\,
        k \in \mathcal{K} \setminus \{N-1\} 
        \label{eq:MVGP_ct_minsep}
    \\&
    y^{(i,j)}_{\Delta r_{\max},k|j} = y^{(i,j)}_{\Delta r_{\max},k+1|j}
        \quad \forall (i,j) \in \mathcal{C},\,
        k \in \mathcal{K} \setminus \{N-1\} 
        \label{eq:MVGP_ct_maxsep}
\end{align}
\end{subequations}
Problem~\eqref{eq:MVGP_ct} resembles~\eqref{eq:MVGP}, with the differences being the dynamics continuity now enforced for the augmented state $\Xbold^{\rm g}$ rather than the concatenated state $\Xbold$, and the introduction of~\eqref{eq:MVGP_ct_minsep} in place of~\eqref{eq:MVGP_minsep}.
Enforcing continuity and periodic boundary conditions on $y^{(i,j)}_{k|j}$ and $y^{(i,j)}_{k+1|j}$ for each pair $(i,j)$ across each consecutive time-steps $k$ to $k+1$ through constraints~\eqref{eq:MVGP_ct_dynamics} and~\eqref{eq:MVGP_ct_minsep} ensure~\eqref{eq:pc_MVGP_integral} is satisfied.
Note that problem~\eqref{eq:MVGP_ct} is still a non-convex NLP due to the dynamics constraints~\eqref{eq:MVGP_ct_dynamics}.

To circumvent linear inequality constraint qualification (LICQ) issues due to linear dependence of gradients of~\eqref{eq:MVGP_ct_dynamics},~\eqref{eq:MVGP_ct_minsep} and~\eqref{eq:MVGP_ct_maxsep}, we adopt the relaxation from Elango et al.~\cite{Elango2024} and replace~\eqref{eq:MVGP_ct_minsep} and~\eqref{eq:MVGP_ct_maxsep} with 
\begin{subequations}    \label{eq:MVGP_ct_minmaxsep_relaxed}
\begin{align}
    y^{(i,j)}_{\Delta r_{\min},k+1|j} - y^{(i,j)}_{\Delta r_{\min},k|j} &\leq \epsilon_{\rm LICQ}
        \quad \forall i, \, j \text{ s.t. } i\neq j, \,
        k \in \mathcal{K} \setminus \{N-1\} 
    \\
    y^{(i,j)}_{\Delta r_{\max},k+1|j} - y^{(i,j)}_{\Delta r_{\max},k|j} &\leq \epsilon_{\rm LICQ}
        \quad \forall i, \, j \text{ s.t. } i\neq j, \,
        k \in \mathcal{K} \setminus \{N-1\} 
    ,
\end{align}
\end{subequations}
where $\epsilon_{\rm LICQ}$ is a user-defined small positive constant.
The chosen value of $\epsilon_{\rm LICQ}$ can be associated with the worst-case violation of the path constraints~\cite{Elango2024}.
Considerations for tuning $\epsilon_{\rm LICQ}$, along with other hyperparameters, are discussed in a subsequent section.

\subsection{Sequential Convex Programming with Augmented Lagrangian}
We adopt an SCP approach to solve the non-convex NLP~\eqref{eq:MVGP_ct}.
Specifically, we adopt the so-called \verb|SCvx*| algorithm~\cite{Oguri2023}, which consists of an Augmented Lagrangian framework to gradually penalize violation of convexified constraints while iteratively adjusting the trust-region on variables.
Having integrated the continuous-time path constraint by augmenting the dynamics and enforcing the relaxed boundary conditions~\eqref{eq:MVGP_ct_minmaxsep_relaxed}, the Augmented Lagrangian framework can readily accommodate the MVGP.

On the $i^{\rm th}$ iteration, we solve a convexified version of~\eqref{eq:MVGP_ct}, where the dynamics constraints~\eqref{eq:MVGP_ct_dynamics} is linearized and relaxed with slack variables $\xibold_{k} \in \mathbb{R}^{6M + q}$ as 
\begin{equation}    \label{eq:linearized_dynamics_contraints}
    \Xbold^{\mathrm{aug}}_{k+1|j} - 
    \left(
        \Phibold^{(i)}(t_{k+1|j},t_{k|j}) \left(
            \Xbold^{\mathrm{aug}}_{k|j} + \begin{bmatrix}
                \Ebold \\ \boldsymbol{0}_{q \times 3M}
            \end{bmatrix} \Ubold_{k|j}
        \right) 
        + \cbold_{k|j}^{(i)}
    \right)
    =
    \xibold_{k}
    \quad 
    \forall k \in \mathcal{K} \setminus \{N-1\} 
\end{equation}
where $\Phibold^{(i)}(t_{k+1|j},t_{k|j})$ is the state-transition matrix of the augmented system obtained by solving the initial value problem
\begin{equation}    
    \begin{cases}
        \dot{\Phibold}^{(i)}(t,t_{k|j}) =
        \left.\dfrac{\partial \Fbold^{\rm g}}{\partial \Xbold^{\rm g}}\right|_{\bar{\Xbold}^{\mathrm{aug},(i)}_{k|j},\bar{\Ubold}^{(i)}_{k|j}}
        \Phibold^{(i)}(t,t_{k|j}),
        \\
        \Phibold^{(i)}(t_{k|j},t_{k|j}) = \Ibold_{6M + q},
    \end{cases}
\end{equation}
and $\cbold^{(i)}_{k|j}$ is given by
\begin{equation}    
    \begin{aligned}
        \cbold_{k|j}^{(i)} &= 
        \bar{\Xbold}^{\mathrm{aug},(i)}_{k+1,k}
        - \Phibold^{(i)}(t_{k+1|j},t_{k|j})
        \left(
        \bar{\Xbold}^{\mathrm{aug},(i)}_{k|j}
        + \begin{bmatrix}
                \Ebold \\ \boldsymbol{0}_{q \times 3M}
            \end{bmatrix} \bar{\Ubold}^{(i)}_{k|j}
        \right)
        ,
        \\
        \bar{\Xbold}^{\mathrm{aug},(i)}_{k+1,k} &= 
        \bar{\Xbold}^{\mathrm{aug},(i)}_{k|j} + \int_{t_{k|j}}^{t_{k+1|j}} \Fbold^{\rm g}(\Xbold^{\rm g},t) \mathrm{d}t
        .
    \end{aligned}
\end{equation}
In addition to the convexified dynamics constraints~\eqref{eq:linearized_dynamics_contraints} and the other convex constraints~\eqref{eq:MVGP_ct_initialconditions} $\sim$ \eqref{eq:MVGP_ct_minsep}, we impose a trust region constraint on the states,
\begin{equation}    \label{eq:trust_region_constraint}
    \| \Xbold^{\rm g}_{k|j} - \bar{\Xbold}^{\mathrm{aug},(i)}_{k|j} \|_{\infty} \leq \Delta_{\rm tr}
    \quad \forall k \in \mathcal{K} \setminus \{N-1\}
    ,
\end{equation}
where $\Delta_{\rm tr} > 0$ is the trust region, thereby preventing the convex problem's solution from lying too far from the reference, where the linearization would yield poor approximations.
Note that the issue of artificial infeasibility due to imposing~\eqref{eq:trust_region_constraint} is circumvented by the introduction of~$\xibold_{k}$, which acts as a unified virtual control/buffer term.
The objective~\eqref{eq:MVGP_ct_obj} is augmented with a penalty function $P(\cdot) \geq 0$, given by
\begin{equation}    
    P(\xibold_{0},\ldots,\xibold_{N-1},\lambdabold_{0}^{(i)},\ldots,\lambdabold_{N-1}^{(i)},w^{(i)}) =
    \sum_{k=1}^{N-1} {\lambdabold_{k}^{(i)}}^T \xibold_{k} + \dfrac{w^{(i)}}{2} \sum_{k=1}^{N-1}\xibold_{k}^T \xibold_{k}
    ,
\end{equation}
where $\lambdabold^{(i)}_{k} \in \mathbb{R}^{6M + q}$ are the Lagrange multipliers on the $k^{\rm th}$ dynamics constraints~\eqref{eq:MVGP_ct_dynamics}, and $w^{(i)} \geq 0$ is a scalar weight on the quadratic term.
Update rules for $\lambdabold^{(i)}_{k}$ and $w^{(i)}$ as well as conditions for accepting and updating the reference solution $\bar{\Xbold}^{\mathrm{aug},(i)}_{k|j}$ and $\bar{\Ubold}^{(i)}_{k|j}$ for $k=1,\ldots,N$ are given in Oguri~\cite{Oguri2023}.


\section{Recursive simulation setup}

We consider a recursive simulation setup similar to our previous work~\cite{Shimane2024-ep} for a single spacecraft guidance problem, but extended to the multi-vehicle case.
At the initial time $t_0$ of the simulation, the true state of the $i^{\rm th}$ spacecraft each spacecraft $\xbold_i(t_0)$ is initialized by shifting by a pre-defined offset $\Delta \xbold_i$ with respect to a baseline NRHO state $\xbold_{\rm bsln} = [\rbold_{\rm bsln}^T, \vbold_{\rm bsln}^T]^T$ and appending an initial insertion error $\delta \xbold_{0} = [\delta \rbold_{0}^T, \delta \vbold_{0}^T]^T$,
\begin{equation}
    \xbold_i(t_0) = \xbold_{\rm bsln} + \Delta \xbold_i + \begin{bmatrix}
        \delta \rbold_{0} \\ \delta \vbold_{0}
    \end{bmatrix}
    ,\quad
    \delta \rbold_{0} \sim \mathcal{N}(\boldsymbol{0}_{3\times1}, \sigma_{r_0}^2\boldsymbol{1}_{3\times1})
    ,\,\,
    \delta \vbold_{0} \sim \mathcal{N}(\boldsymbol{0}_{3\times1}, \sigma_{v_0}^2\boldsymbol{1}_{3\times1})
    ,
\end{equation}
where $\delta \rbold_{0}$ and $\delta \vbold_{0}$ are realized independently for each $i$.
The initial offset $\Delta \xbold_i$ is defined with a fixed separation; specifically, for the two-spacecraft case,
\begin{equation}
    \Delta \xbold_0 = \begin{bmatrix}
        0 & 0 & 1.5 \Delta r_{\min}/2 & 0 & 0 & 0 
    \end{bmatrix}^T
    ,\quad
    \Delta \xbold_1 = \begin{bmatrix}
        0 & 0 & -1.5 \Delta r_{\min}/2 & 0 & 0 & 0 
    \end{bmatrix}^T
    .
\end{equation}
Each time the MPC is invoked, the state estimate for the $i^{\rm th}$ spacecraft, $\hat{\xbold}_i(t_j)$, is approximated by appending a navigation error $\delta \xbold_{\rm nav} = [\delta \rbold_{\rm nav}^T, \delta \vbold_{\rm nav}^T]^T$,
\begin{equation}
    \hat{\xbold}_i(t_j)
    = 
    \xbold_i(t_j) +
    \begin{bmatrix}
        \delta \rbold_{\rm nav} \\ \delta \vbold_{\rm nav}
    \end{bmatrix}
    ,\quad
    \delta \rbold_{\rm nav} \sim \mathcal{N}(\boldsymbol{0}_{3\times1}, \sigma_{\rm nav}^2\boldsymbol{1}_{3\times1})
    ,\,\,
    \delta \vbold_{\rm nav} \sim \mathcal{N}(\boldsymbol{0}_{3\times1}, \sigma_{\rm nav}^2\boldsymbol{1}_{3\times1})
    ,
\end{equation}
where $\delta \rbold_{\rm nav}$ and $\delta \vbold_{\rm nav}$ are realized independently for each $i$.
The control maneuver computed by the MPC for the $i^{\rm th}$ spacecraft, $\ubold_{i,0|j}$, is also subject to errors. Let $\tilde{\ubold}_{i,0|j}$ denote the maneuver physically imparted on the spacecraft. 
Following Gates' model~\cite{Gates1963}, we realize independently directional, relative magnitude, and absolute magnitude errors on $\ubold_{i,0|j}$, such that
\begin{equation}
\begin{aligned}
    \tilde{\ubold}_{i,0|j} &= \Tbold(\delta \phi) \left( \ubold_{i,0|j}\left(1+\delta u_{\rm rel} \right) + \delta \ubold_{\rm abs} \right)
    ,\,\,\,
    \delta u_{\rm rel} \sim \mathcal{N}(0,\sigma_{\rm rel}^2), \,\,
    \delta \ubold_{\rm abs} \sim \mathcal{N}(\boldsymbol{0}_{3\times1},\sigma_{\rm abs}^2\boldsymbol{1}_{3\times1})
    \\
    \Tbold(\delta \phi) &= \cos (\delta \phi) \boldsymbol{I}_3+\sin (\delta \phi) \boldsymbol{i}^{\times} +[1-\cos (\delta \phi)] \boldsymbol{i} \boldsymbol{i}^T, \quad \delta \phi \sim \mathcal{N}\left(0, \sigma_{\varphi}^2\right),
\end{aligned}
\end{equation}
where $\boldsymbol{i}$ is a random unit vector, and $\boldsymbol{i}^{\times}$ is its skew-symmetric form.

\section{Numerical Results}
We compare guidance results (a) without any path constraints, (b) with discrete-time path constraints enforced at the control nodes, and (c) with path constraints enforced continuously in time.
Case (a) involves solving the NLP~\eqref{eq:MVGP} without constraints~\eqref{eq:MVGP_minsep}, case (b) involves replacing constraints~\eqref{eq:MVGP_minsep} with~\eqref{eq:minsep_nodes}, and case (c) involves solving~\eqref{eq:MVGP_ct}.
The dynamics and SCP framework are implemented in Julia.
The convex program within the SCP is solved using Clarabel~\cite{Goulart2024}.
Parameters for the MVGP, the corresponding MPC, along with hyperparameters for the SCP scheme, are given in Table~\ref{tab:mvgp_parameters}.

\begin{table}[t]
\centering
\caption{MVGP parameters and SCP hyperparameters}
\label{tab:mvgp_parameters}
\begin{tabular}{@{}ll@{}}
    \toprule
    Parameter & Value \\ \midrule
    Control true anomaly $\theta_{\rm man}$, \SI{}{deg} & $180^{\circ}$ \\
    MPC prediction horizon $N$ & $6$ \\
    Targeting ellipsoid position radius $\epsilon_r$, \SI{}{km}  & cooperative: $20$ / hierarchical: $30$ \\
    Targeting ellipsoid velocity radius $\epsilon_v$, \SI{}{m/s} & $5$ \\
    Minimum separation radius $\Delta r_{\min}$, \SI{}{km}       & $10$ \\
    Minimum separation radius buffer $\delta r_{\min}$ \SI{}{km} & $10$ \\
    Maximum separation radius $\Delta r_{\max}$, \SI{}{km}       & $200$ \\
    Maximum separation radius buffer $\delta r_{\max}$ \SI{}{km} & $25$ \\
    SCP initial weight $w^{(0)}$ & $10^2$ \\
    Continuous-time constraint KKT relaxation parameter $\epsilon_{\rm LICQ}$ & $10^{-6}$ \\
    Path constraint reformulation parameter $\kappa$ & $10^5$ \\
    \bottomrule
\end{tabular}
\end{table}
 
\begin{table}[t]
\centering
\caption{Error parameters for recursive simulation}
\label{tab:mvgp_parameters}
\begin{tabular}{@{}ll@{}}
    \toprule
    Parameter & Value \\ \midrule
    Initial position error $3\text{-}\sigma_{r_0}$, \SI{}{km}   & $1.0$ \\ 
    Initial velocity error $3\text{-}\sigma_{v_0}$, \SI{}{cm/s} & $1.0$ \\ 
    Navigation position error $3\text{-}\sigma_{r,\mathrm{nav}}$, \SI{}{km}    & $0.1$ \\ 
    Navigation position error $3\text{-}\sigma_{r,\mathrm{nav}}$, \SI{}{cm/s}  & $1.0$ \\ 
    Maneuver relative magnitude error $3\text{-}\sigma_{u,\mathrm{rel}}$, $\%$ & $1.5$ \\
    Maneuver absolute magnitude error $3\text{-}\sigma_{u,\mathrm{abs}}$, \SI{}{mm/s} & $1.42$ \\
    Maneuver execution direction error $3\text{-}\sigma_{u,\mathrm{dir}}$, \SI{}{deg} & $0.5^{\circ}$ \\
    \bottomrule
\end{tabular}
\end{table}

\subsection{Recursive Guidance under Uncertainty}
\begin{figure}[t]
     \centering
     \begin{subfigure}[b]{0.48\textwidth}
         \centering
         \includegraphics[width=\textwidth]{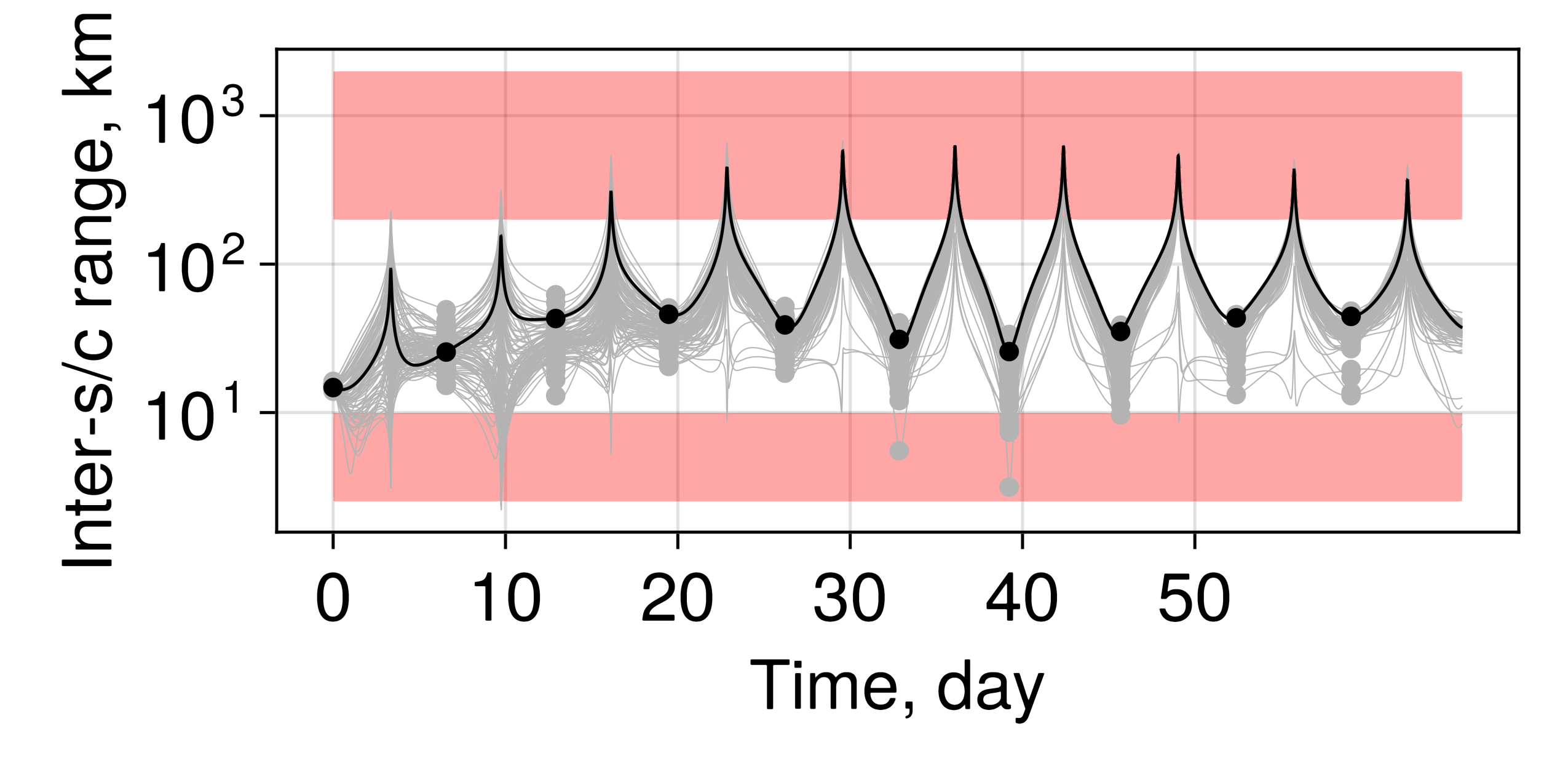}
         \caption{No path constraints}
         \label{fig:cooperative_problem_range_history_case0}
     \end{subfigure}
     \hfill
     \begin{subfigure}[b]{0.48\textwidth}
         \centering
         \includegraphics[width=\textwidth]{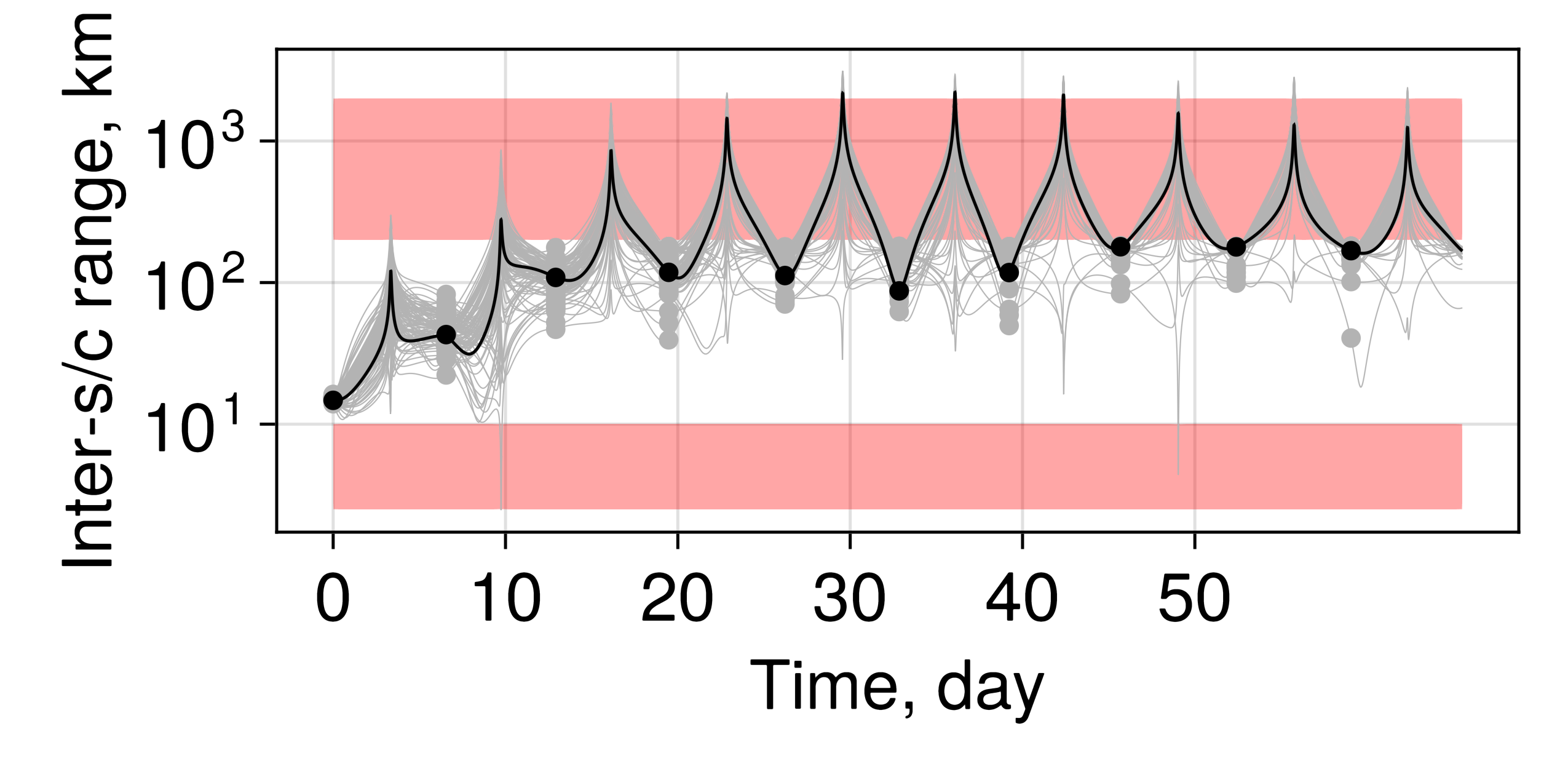}
         \caption{Path constraints at control nodes}
         \label{fig:cooperative_problem_range_history_case1}
     \end{subfigure}
     \\
     \begin{subfigure}[b]{0.48\textwidth}
         \centering
         \includegraphics[width=\textwidth]{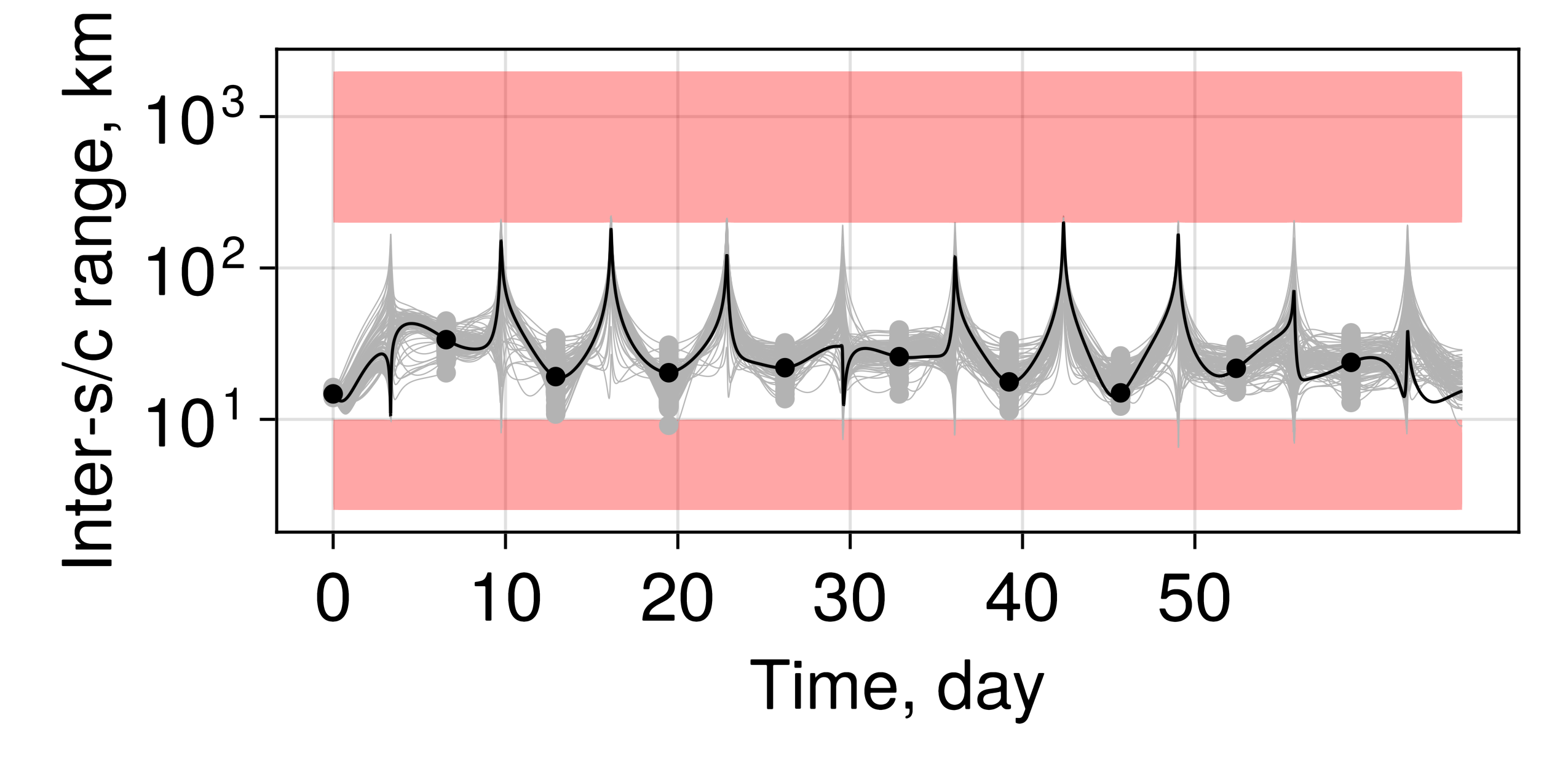}
         \caption{Continuous-time path constraints}
         \label{fig:cooperative_problem_range_history_case2}
     \end{subfigure}
    \caption{Monte Carlo samples of inter-spacecraft range recursively controlled with MPC}
    \label{fig:cooperative_problem_range_history}
\end{figure}

We consider a Monte-Carlo experiment incorporating initial injection, navigation, and control execution errors for the formation flight of two spacecraft for a duration of 10 revolutions along the nominal NRHO.
Figure~\ref{fig:cooperative_problem_range_history} shows the inter-spacecraft range for 100 Monte-Carlo samples with (a), (b), and (c).
In the Figure, the circle markers denote the control nodes.
With (a) no path constraint, the formation violates the separation constraints nearly every revolution; in contrast, enforcing the separation constraints (b) at control nodes or (c) continuously significantly reduces instances of violation.
Due to the presence of uncertainties, some of the samples from (b) and (c) still violate the minimum separation, albeit to a much smaller extent and with lower likelihood.
Comparing cases (b) and (c), the range history reveals that the recursive solution traces different trends; with (b), the relative trajectory follows a path such that the maximum separation constraint is active at the nodes; with (c), with the trajectory bounded by both the minimum and maximum separation at all times, the range is closer to the minimum separation limit at the nodes.

\begin{figure}[t]
     \centering
     \begin{subfigure}[b]{0.48\textwidth}
         \centering
         \includegraphics[width=\textwidth]{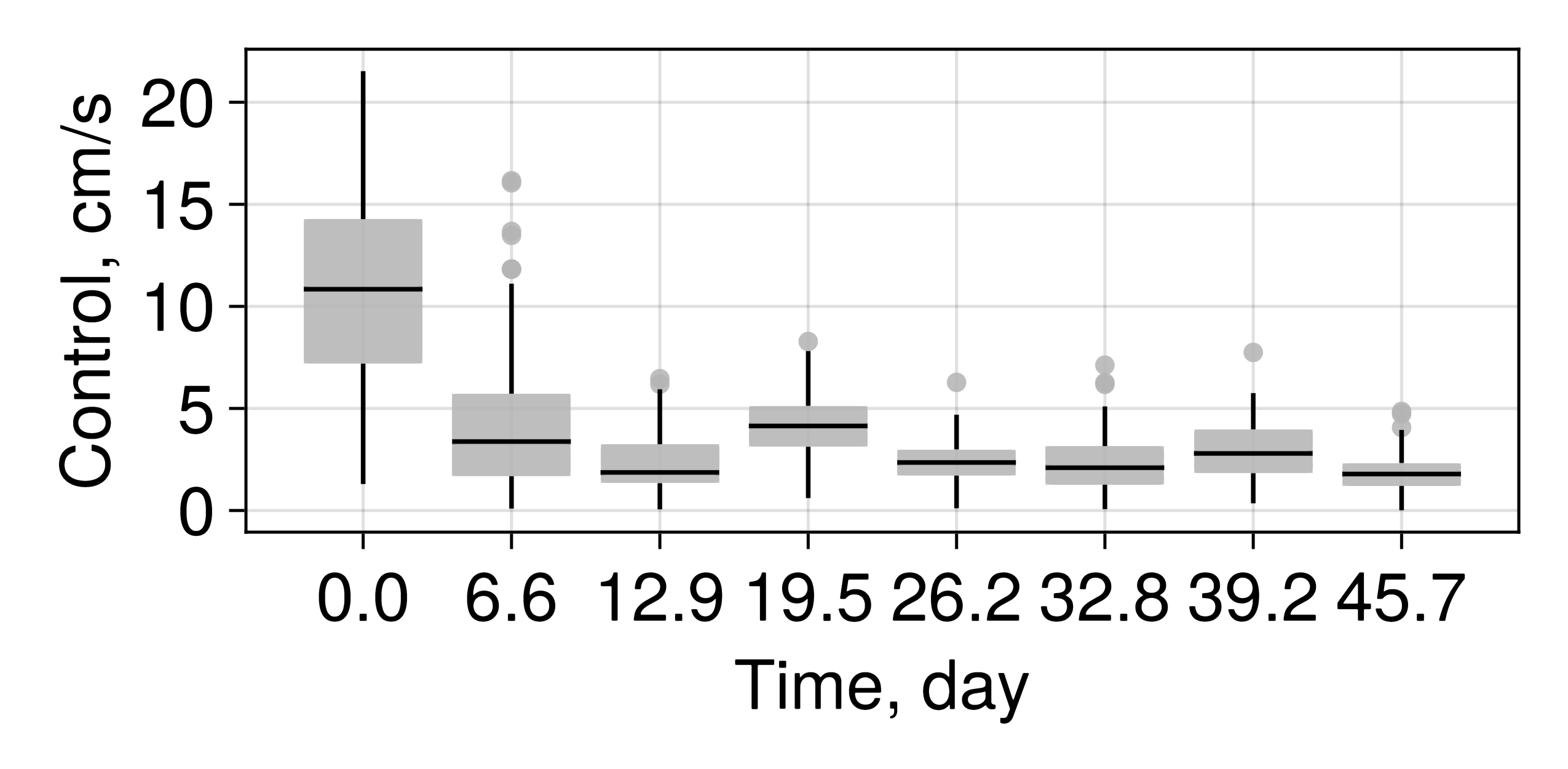}
         \caption{No path constraints}
         \label{fig:cooperative_problem_umag_case0}
     \end{subfigure}
     \hfill
     \begin{subfigure}[b]{0.48\textwidth}
         \centering
         \includegraphics[width=\textwidth]{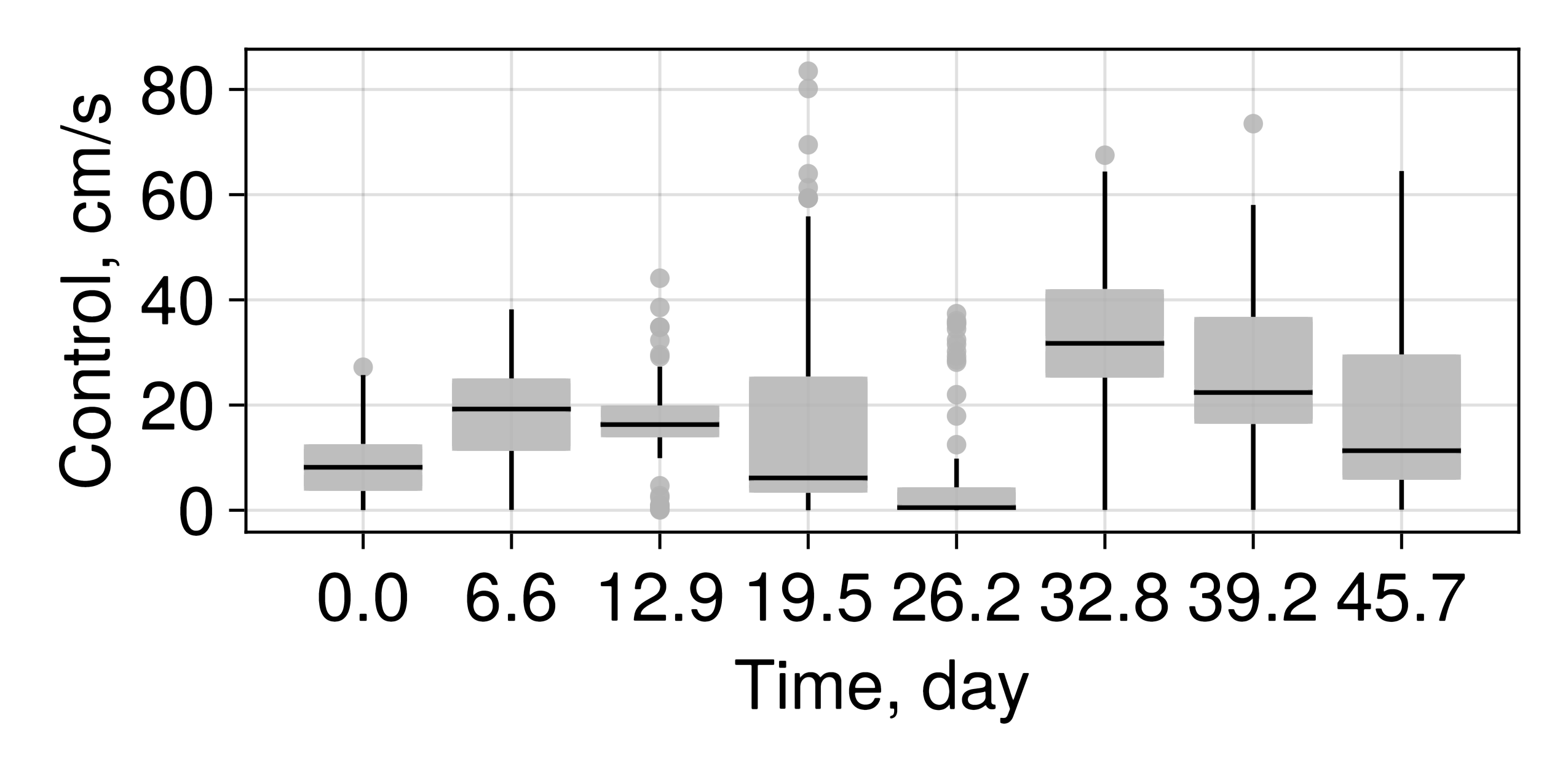}
         \caption{Path constraints at control nodes}
         \label{fig:cooperative_problem_umag_case1}
     \end{subfigure}
     \\
     \begin{subfigure}[b]{0.48\textwidth}
         \centering
         \includegraphics[width=\textwidth]{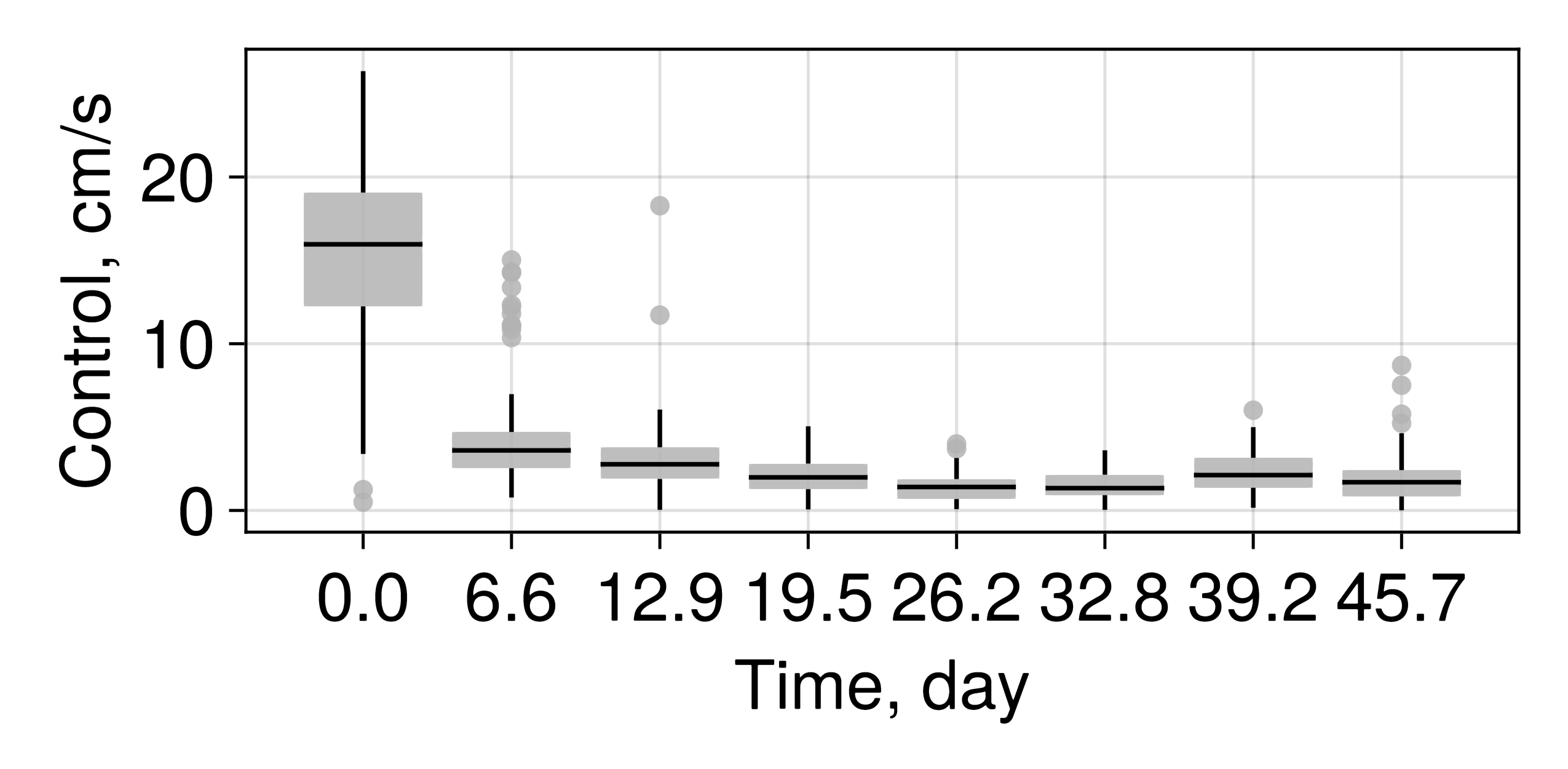}
         \caption{Continuous-time path constraints}
         \label{fig:cooperative_problem_umag_case2}
     \end{subfigure}
    \caption{Monte Carlo samples of executed control magnitudes}
    \label{fig:cooperative_problem_umag}
\end{figure}

The distribution of the executed control from the Monte Carlo samples is shown in Figure~\ref{fig:cooperative_problem_umag}.
While the distribution history for (a) and (c) is similar, the cost distribution for (b) has notable differences.
With (a) and (c), besides the initial control magnitude, which has a relatively large distribution, the control magnitude in the second recursion onward is much smaller.
The first maneuver is impacted by the realization of the initial insertion error $\delta \xbold_0$, thus resulting in a larger distribution; meanwhile, the subsequent maneuver magnitudes have much smaller distributions, suggesting that the controlled trajectory is in a steady-state regime.
The low control distribution with (c) is notable, as it suggests that the MPC has entered a cheaply controllable configuration under uncertainty that abides by both the minimum and maximum separation constraints.
In contrast, (b) results in maneuver magnitudes having larger distributions on a number of recursions, suggesting that the trajectories are more susceptible to uncertainties.

\begin{figure}[t]
     \centering
     \begin{subfigure}[b]{0.48\textwidth}
         \centering
         \includegraphics[width=\textwidth]{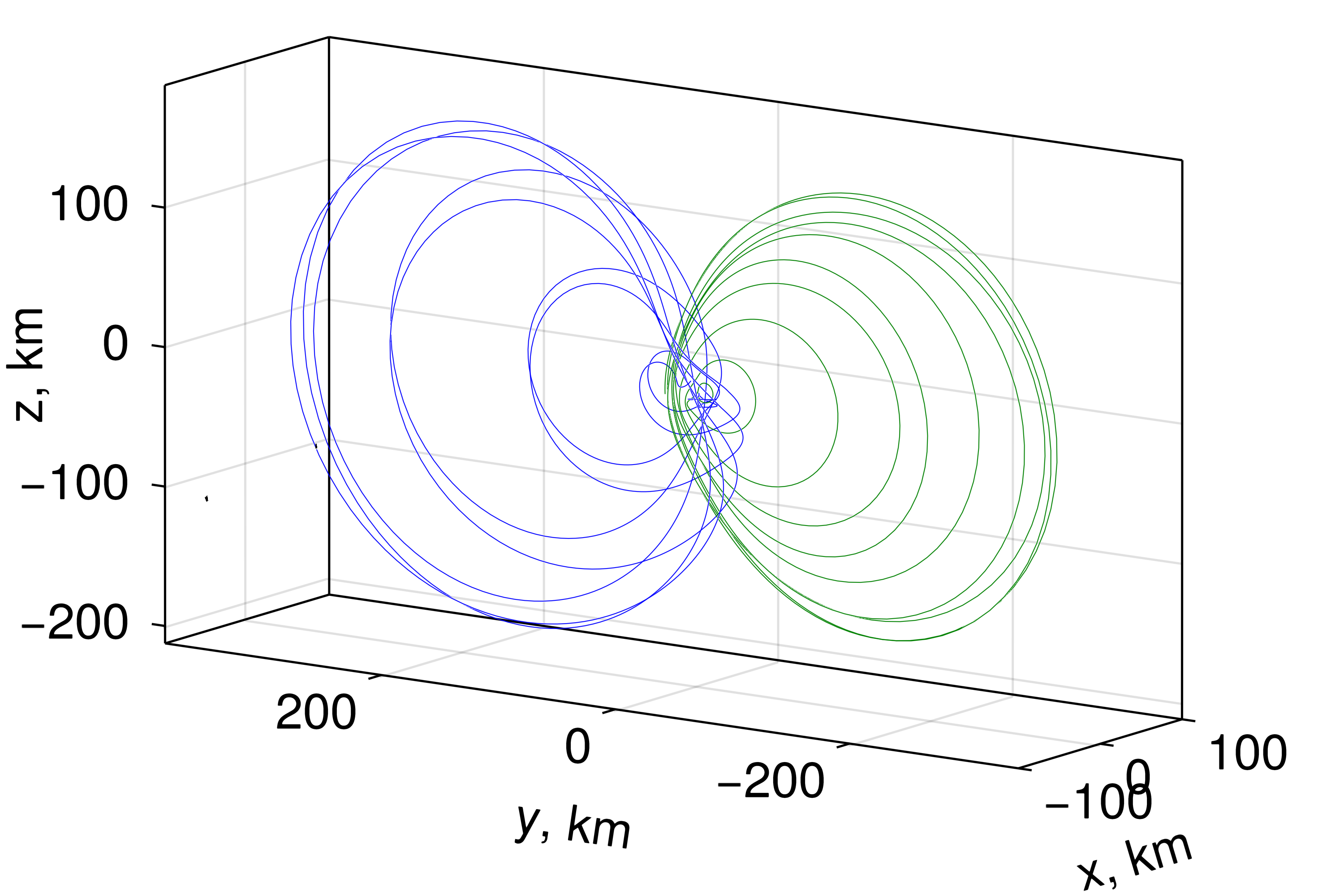}
         \caption{No path constraints}
         \label{fig:cooperative_problem_reltraj_case0}
     \end{subfigure}
     \hfill
     \begin{subfigure}[b]{0.48\textwidth}
         \centering
         \includegraphics[width=\textwidth]{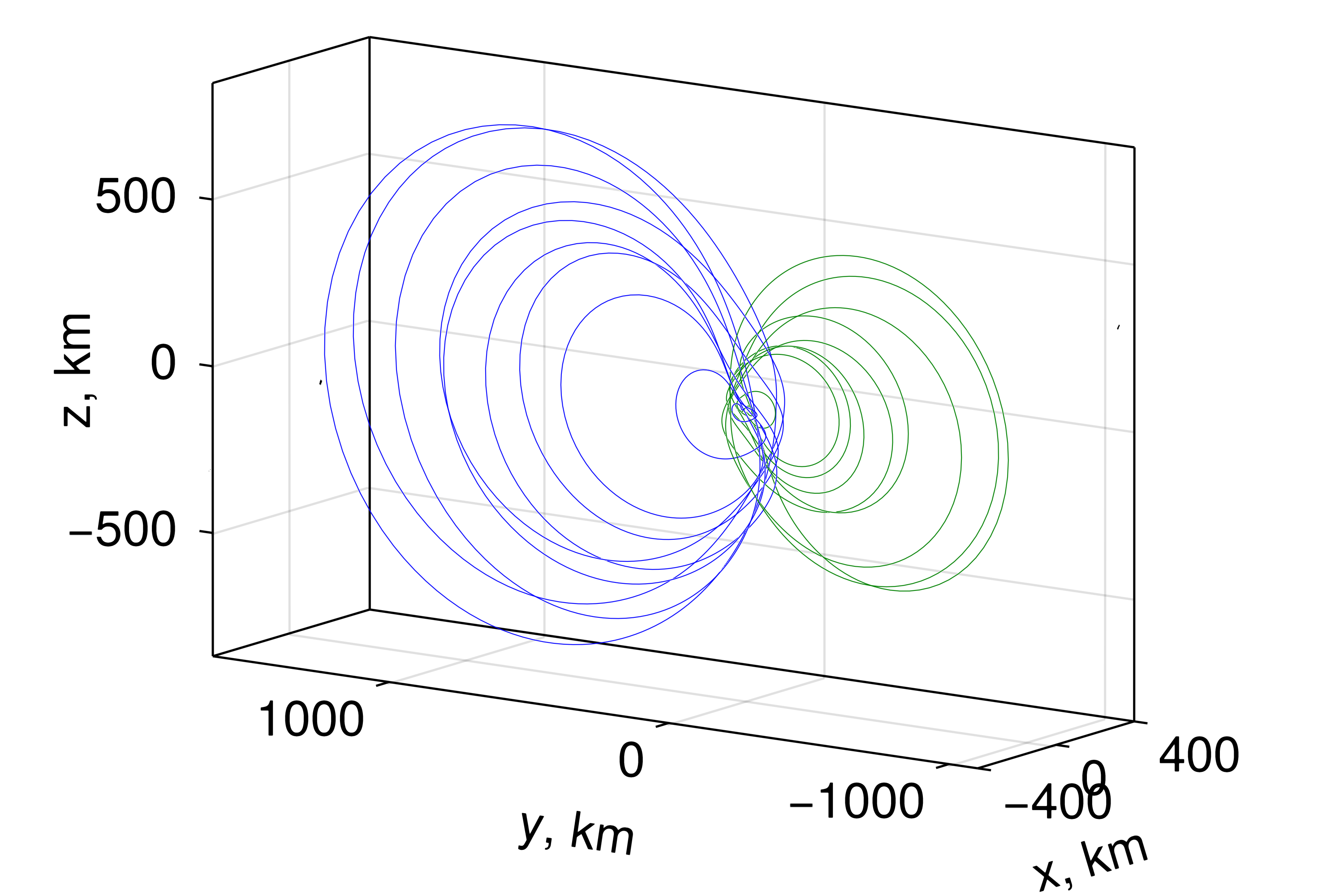}
         \caption{Path constraints at control nodes}
         \label{fig:cooperative_problem_reltraj_case1}
     \end{subfigure}
     \\
     \begin{subfigure}[b]{0.48\textwidth}
         \centering
         \includegraphics[width=\textwidth]{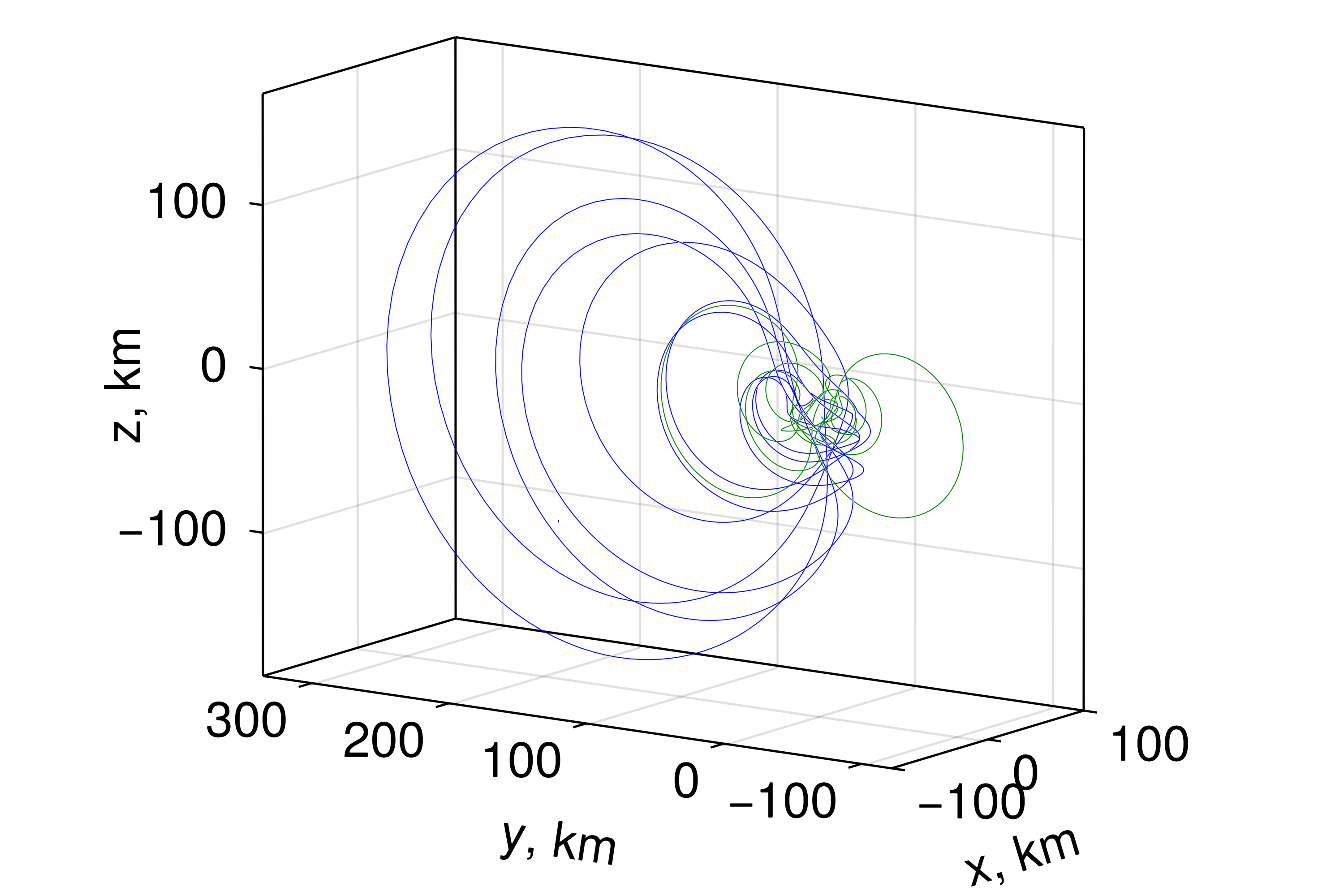}
         \caption{Continuous-time path constraints}
         \label{fig:cooperative_problem_reltraj_case2}
     \end{subfigure}
    \caption{Relative trajectory over a single Monte Carlo sample in Earth-Moon rotating frame}
    \label{fig:cooperative_problem_reltraj}
\end{figure}

Figure~\ref{fig:cooperative_problem_reltraj} shows a sample relative trajectory for the three cases.
In all cases, the relative trajectories suggest that the path is steered along a quasi-periodic torus, even though no explicit information is provided to the controller.
As informed by the inter-spacecraft range in Figure~\ref{fig:cooperative_problem_range_history_case1} and illustrated in Figure~\ref{fig:cooperative_problem_reltraj_case1}, the relative trajectory with case (b) follows a much larger amplitude compared to the other two cases.
Also notable is case (c), where the amplitude of the first and second spacecraft is substantially different. This strategy allows for the continuously bounded inter-spacecraft range, as shown in Figure~\ref{fig:cooperative_problem_range_history_case2}.

\section{Conclusions}
This work addressed a guidance scheme for formation flying on LPOs.
We developed an MPC where a non-convex optimal control problem with continuous-time path constraints on the inter-spacecraft range is recursively solved via SCP.
The bounds on the path constraints are defined with a tightening factor across the MPC's control horizon to empirically ensure the MPC yields control sequences that are recursively feasible under uncertainty.
We successfully demonstrated the proposed MPC with Monte Carlo simulations along the 9:2 resonant NRHO in the HFEM with one control maneuver per revolution, incorporating realistic uncertainty models.
We find that the use of continuous-time path constraints pushes the formation towards a configuration that yields cheaper cumulative control costs compared to enforcing the path constraints only at the control nodes.




\bibliographystyle{AAS_publication}   
\bibliography{references}   

\end{document}